\documentclass[12pt,reqno]{amsart} 
\usepackage{amscd,amssymb,amsmath}
\usepackage[arrow,matrix,curve]{xy}

\newtheorem{theorem}{Theorem}[section]
\newtheorem{proposition}[theorem]{Proposition}
\newtheorem{lemma}[theorem]{Lemma}
\newtheorem{corollary}[theorem]{Corollary}
\theoremstyle{definition}
\newtheorem{definition}[theorem]{Definition}

\newtheorem{remark}[theorem]{Remark}

\numberwithin{equation}{section} 
\numberwithin{figure}{section}
\numberwithin{table}{section}

\DeclareMathOperator{\syscat}{{\rm cat_{sys}}}
\DeclareMathOperator{\Indet}{{\rm Indet}}

\def\N {{\mathbb N}}

\def\R {{\mathbb R}}

\def\Z {{\mathbb Z}}

\def\gmetric {{\mathcal G}}

 \def\Hom{\operatorname{Hom}}

\def\rk{\operatorname{{rank}}} \def\dR{{\operatorname{dR}}}

\def\IQ{\operatorname{IQ}}

\newcommand{\integerlat}[1]{{H^{#1}_{\dR}(X,\Z)}}

\def\sys{{\rm sys}}  

\DeclareMathOperator{\stsys}{{\rm stsys}}

\def\vol{{\rm vol}}    \def\length{{\rm length}}

\def\etc {{\it etc}} 

\def\cf {{\it cf.\ }} \def\ie {{\it i.e.\ }} 

\def\vs {{\it vs.\ }}

\def\two {{m}} \def\five {{3m-1}} \def\three {{2m}}
\def\four {{2m}}
\def\x {{x_0}}

\begin{document}

\title [Systolic inequalities and Massey products] {Systolic
inequalities and Massey products in simply-connected manifolds
}

\author[M.~Katz]{Mikhail Katz$^*$} \address{Department of Mathematics,
Bar Ilan University, Ramat Gan 52900 Israel}
\email{katzmik@math.biu.ac.il} \thanks{$^*$Supported by the Israel
Science Foundation (grants no.\ 620/00, 84/03, and 1294/06)}


\begin{abstract}
We show that the existence of a nontrivial Massey product in the
cohomology ring~$H^*(X)$ imposes global constraints upon the
Riemannian geometry of a manifold~$X$.  Namely, we exhibit a suitable
systolic inequality, associated to such a product.  This generalizes
an inequality proved in collaboration with Y.~Rudyak, in the case
when~$X$ has unit Betti numbers, and realizes the next step in
M.~Gromov's program for obtaining geometric inequalities associated
with nontrivial Massey products.  The inequality is a volume lower
bound, and depends on the metric via a suitable isoperimetric
quotient.  The proof relies upon W.~Banaszczyk's upper bound for the
successive minima of a pair of dual lattices.  Such an upper bound is
applied to the integral lattices in homology and cohomology of~$X$.
The possibility of applying such upper bounds to obtain volume lower
bounds was first exploited in joint work with V.~Bangert.  The latter
work deduced systolic inequalities from nontrivial {\em cup-product\/}
relations, whose role here is played by Massey products.
\end{abstract}



\maketitle
\tableofcontents

\section{Volume bounds and systolic category}

A general framework for systolic geometry in a topological context was
proposed in \cite{KR1}, in terms of a new invariant called {\em
systolic category}, denoted~$\syscat(X)$, of a space~$X$.  The
terminology is inspired by the intriguing connection which emerges
with the classical invariant called the Lusternik-Schnirelmann
category.  Namely, the two categories (\ie the two integers) coincide
for 2-complexes \cite{KRS}, as well as for 3-manifolds, orientable or
not \cite{KR1, KR2}, attain their maximal value simultaneously, both
admit a lower bound in terms of real cup-length, both are sensitive to
Massey products, \etc.

\begin{definition}
The stable~$k$-systole of a Riemannian manifold is the least stable
norm of a nonzero element in the integer lattice in its
$k$-dimensional homology group with real coefficients.
\end{definition}
A more detailed definition appears below, \cf formula \eqref{34}.

The invariant~$\syscat$ is defined in terms of the existence of volume
lower bounds of a certain type.  Namely, these are bounds by products
of lower-dimensional systoles.  The invariant~$\syscat$ is, roughly,
the greatest length~$d$ of a product
\[
\prod_{i=1}^d \sys_{k_i}
\]
of systoles which provides a universal lower bound for the volume, \ie
a curvature-independent lower bound of the following form:
\[
\prod_{i=1}^d \sys_{k_i}(\gmetric) \leq C \vol(\gmetric),
\]
see~\cite{KR1} for details.  The definitions of the systolic
invariants involved may also be found in \cite{Gr1, CK, KL}.

We study stable systolic inequalities satisfied by an arbitrary
metric~$\gmetric$ on a closed, smooth manifold~$X$.  We aim to go
beyond the multiplicative structure, defined by the cup product, in
the cohomology ring, whose systolic effects were studied in
\cite{Gr1, He, BK1, BK2}, and explore the systolic influence of
Massey products.

\begin{remark}
This line of investigation is inspired by M. Gromov's remarks
\cite[7.4.$C'$, p.~96]{Gr1} and \cite[7.5.C, p.~102]{Gr1}, outlining a
program for obtaining geometric inequalities associated to nontrivial
Massey products of any length.  The first step in the program was
carried out in \cite{KR1} in the presence of a nontrivial triple
Massey product in a manifold with unit Betti numbers.
\end{remark}

In the present work, we exploit W.~Banaszczyk's bound \eqref{szcz} for
the successive minima of a pair of dual lattices, applied to the
integral lattices in homology and cohomology of~$X$.  The possibility
of exploiting such bounds to obtain inequalities was first
demonstrated in joint work with V. Bangert \cite{BK1} on systolic
inequalities associated to nontrivial cup product relations in the
cohomology ring of~$X$.

Whenever a manifold admits a nontrivial Massey product, we seek to
exhibit a corresponding inequality for the stable systoles.  While
nontrivial cup product relations in cohomology entail stable systolic
inequalities which are metric-independent and curvature-free
\cite{BK1}, the influence of Massey products on systoles is more
difficult to pin down.  The inequalities obtained so far do depend
mildly on the metric, via isoperimetric quotients, \cf \eqref{isop}.

The idea is to show that if, in a certain dimension~$k \leq n$, one
can span the cohomology by classes which can be expressed in terms of
lower-dimensional classes by either Massey or cup products, then the
stable~$k$-systole (\cf Definition~\ref{stsys}) admits a bound from
below in terms of lower-dimensional stable systoles, and of certain
isoperimetric constants of the metric, but no further metric data.
Typical examples are inequalities \eqref{sys2a}, \eqref{sys2},
\eqref{82}.

Massey products and isoperimetric quotients are reviewed in
Section~\ref{two}.  The theorems are stated in Section~\ref{thms}.
Banaszczyk's results are reviewed in Section~\ref{four}.  The key
notion of quasiorthogonal element of a Massey product is defined in
Section~\ref{sect:proofs}.  The theorems are proved in
Section~\ref{five}.

The basic reference for this material is M. Gromov's monograph
\cite{Gr3}, with additional details in the earlier works~\cite{Gr1,
Gr2}.  For a survey of progress in systolic geometry up to 2003, see
\cite{CK}.  More recent results include a study of optimal
inequalities of Loewner type \cite{Am, IK, BCIK2, KL, KS1}, as well as
near-optimal asymptotic bounds \cite{BB, Ka3, KS2, KS3, Sa1, Sa3,
KSV}, while generalisations of Pu's inequality are studied in
\cite{BCIK1} and \cite{e7}.  For an overview of systolic questions,
see \cite{SGT}.

\section{Massey products and isoperimetric quotients}
\label{two}

In Theorem~\ref{22}, we will use a hypothesis which in the case of no
indeterminacy of Massey products, amounts simply to requiring every
cohomology class to be a sum of Massey products.  In general, the
condition is slightly stronger, and informally can be described by
saying that any system of representatives of Massey products already
spans the entire cohomology space.

Following the notation of \cite{KR1}, consider (homogeneous) cohomology
classes~$u,v, w$ with~$uv=0=vw$.  Then the triple Massey product
\[
\langle u,v,w\rangle \subset H^*_\dR
\]
is defined as follows.  Let~$a,b,c$ be closed differential forms whose
homology classes are~$u,v,w$ respectively. Then~$dx=ab$,~$dy=bc$ for
suitable differential forms~$x,y$.  Then~$\langle u,v,w\rangle$ is
defined to be the set of elements of the form
\[
xc - (-1)^{|u|}ay,
\]
see \cite{M, RT} for more details.  The set~$\langle u,v,w\rangle$ is
a coset with respect to the so-called {\it indeterminacy subgroup}
$\Indet \subset H^{|u|+|v|+|w|-1}$, defined as follows:
\begin{equation}
\label{21}
\Indet = uH^{|v|+|w|-1}+H^{|u|+|v|-1}w.
\end{equation}
A Massey product is said to be {\em nontrivial\/} if it does not
contain~$0$.

\begin{definition}
\label{mt}
Let~$m\geq 1$.  The~$(3m-1)$-dimensional de Rham cohomology space of a
manifold~$X$ is {\em of Massey type\/} if it has the following
property.  Let~$V\subset H^{3m-1}_{\dR}(X)$ be a subspace with
nonempty intersection with every nontrivial triple Massey product
$\langle u,v,w\rangle$,~$u,v,w\in H^m_\dR(X)$.  Then~$V =
H^{3m-1}_{\dR}(X)$.
\end{definition}

Given a compact Riemannian manifold~$(X,\gmetric)$, and a positive
integer~$k\leq \dim X$, denote by~$\IQ_k=\IQ_k(\gmetric)$ the
isoperimetric quotient, defined by
\begin{equation}
\label{isop}
\IQ_k(\gmetric) = \sup_{\alpha \in \Omega^k(X)} \inf_\beta \left\{
\left.  \tfrac{\|\beta \|^*}{ \| \alpha \|^* } \;\right| d\beta =
\alpha \right\},
\end{equation}
where~$\|\;\|^*$ is the comass norm \cite{Fe2}, and the supremum is
taken over all exact~$k$-forms.  The relation of such quotients to
filling inequalities is described in \cite[Section 4, Proposition
1]{Sik}, \cf \cite[item~4.13]{Fe2}.

\section{The results}
\label{thms}

The following theorem generalizes \cite[Theorem~13.1]{KR1} to the case
of arbitrary Betti number.

\begin{theorem}
\label{22}
Let~$X$ be a connected closed orientable smooth manifold.  Let~$m\ge
1$, and assume~$b=b_m(X)>0$.  Furthermore, assume that the following
three hypotheses are satisfied:
\begin{enumerate}
\item 
\label{h1a} 
the cup product map~$\cup\colon H^m_\dR(X) \otimes H^m_\dR (X) \to
H^{2m}_\dR(X)$ is the zero map;
\item 
\label{h2a} 
the space~$H^{3m-1}_{\dR}(X)$ is of Massey type in the sense of
Definition~\ref{mt};
\item
the group~$H^{2m}(X,\Z)$ is torsionfree.
\end{enumerate}
Then every metric~$\gmetric$ on~$X$ satisfies the inequality
\begin{equation}
\label{sys2a}
\stsys_\two(\gmetric)^3 \leq C(m) (b (1+\log b))^3 \IQ _{\three}
(\gmetric) \stsys_\five(\gmetric),
\end{equation}
where~$C(m)$ is a constant depending only on~$m$.
\end{theorem}

Note that the dimensionality of the factor~$\IQ_k(\gmetric)$ is
$(\length)^{+1}$, making inequality \eqref{sys2a} scale-invariant, \cf
\cite [$7.4.C'$, p.~96 and~$7.5.C$, p.~102]{Gr1}.

The proof of Theorem~\ref{22} appears in Section~\ref{five}.

An important special case is a lower bound for the total volume.
While Hypothesis \ref{h2a} of Theorem \ref{22} is rather restrictive,
similar inequalities can be proved in the presence of a nontrivial
Massey product, even if Hypothesis \ref{h2a} is not satisfied,
provided one replaces the systole in the right hand side by the total
volume.  The simplest example of a theorem along these lines is the
following.

\begin{theorem}
\label{222}
Let~$X$ be a closed orientable smooth manifold of dimension~$7$.
Assume that the following three hypotheses are satisfied:
\begin{enumerate}
\item \label{h1} the cup product vanishes on~$H^2_{\dR}(X)$;
\item \label{h3} there are classes~$u,v,w\in H^2_{\dR}(X)$ such that
the triple Massey product~$\langle u,v,w \rangle \subset
H^{5}_{\dR}(X)$ is nontrivial;
\item
the group~$H^4(X,\Z)$ is torsionfree.
\end{enumerate}
Then every metric~$\gmetric$ on~$X$ satisfies the inequality
\begin{equation}
\label{sys2}
\stsys_2(\gmetric)^4 \leq C(b_2(X))\, \IQ_{4}(\gmetric)
\vol_7(\gmetric),
\end{equation}
where the constant~$C(b_2(X))>0$ depends only on the second Betti
number of~$X$.
\end{theorem}
Examples of manifolds to which Theorem~\ref{22} and Theorem~\ref{222}
can be applied, were constructed by A.~Dranishnikov and Y.~Rudyak
\cite{DR}.

Our Theorem~\ref{222} implies the following bound for the
$\IQ$-modified systolic category, \cf \cite[Remark~13.1]{KR1}.

\begin{corollary}
Under the hypotheses of Theorem~\ref{222}, the manifold~$X$ satisfies
the bound~${\rm cat}_{\sys}^{\IQ}(X) \geq 3$.
\end{corollary}

\begin{corollary}
Suppose in addition to the hypotheses of Theorem~\ref{222} that~$X$ is
simply connected.  Then~${\rm cat}_{\sys}^{\IQ}(X) \geq {\rm
cat}_{{\rm LS}}(X)$.
\end{corollary}

\begin{proof}
By \cite[Theorem~1.50]{CLOT}, the Lusternik-Schnirelmann category of
$X$ equals 3.
\end{proof}

Our last result attempts to go beyond both Theorem~\ref{22} and
Theorem~\ref{222}, in the sense of obtaining a lower bound for a
$k$-systole other than the total volume, in a situation where Massey
products do not necessarily span~$k$-dimensional cohomology.

\begin{proposition}
\label{81}
Consider a closed manifold~$X$ with a nontrivial triple Massey product
containing an element~$u\in H^5(X)$.  Assume that the following three
hypotheses are satisfied:
\begin{enumerate}
\item
the cup product vanishes on~$H^2(X)$;
\item
the~$8$-dimensional cohomology of~$X$ is spanned by classes of type~$u
\cup v$ and~$w$, where~$v \in H^3(X)$, while~$w\in H^8(X)$ is the cup
square of a~$4$-dimensional class;
\item
the group~$H^4(X,\Z)$ is torsionfree.
\end{enumerate}
Then every metric~$\gmetric$ on~$X$ satisfies the inequality
\begin{equation}
\label{82}
\min \left\{ \frac {\stsys_2(\gmetric)^3\stsys_3(\gmetric)}
{\IQ_4(\gmetric)}, \; \stsys_4(\gmetric)^2 \right\} \leq C(X)
\stsys_8(\gmetric),
\end{equation}
where~$C(X)>0$ is a constant depending only on the homotopy type
of~$X$.
\end{proposition}

The proof appears in Section~\ref{five}.

\section{Banaszczyk's bound for the successive minima of a lattice}
\label{four}

Let~$B$ be a finite-dimensional Banach space, equipped with a
norm~$\|\;\|$.  Let~$L\subset B$ be a lattice of maximal rank
$\rk(L)=\dim(B)$.  Let~$b=\rk(L)=\dim(B)$.

\begin{definition}
\label{mindef}
For each~$k=1,2,\dots, b$, define the {\em~$k$-th successive
minimum}~$\lambda_k$ of the lattice~$L$ by setting
\begin{equation}
\label{success}
\lambda_k(L,\|\;\|) = \inf\left\{\lambda\in\R \ \left| \,
\begin{array}{l}%
\exists \text{ lin. indep. } v_1, \ldots, v_k \in L \\ \text{\ with\ }
\|v_i\|\leq \lambda,\quad i=1,\ldots,k
\end{array}  
\right. \right\}.
\end{equation}
\end{definition}

In particular, the ``first'' successive minimum,~$\lambda_1
(L,\|\;\|)$, is the least length of a nonzero element in~$L$.

\begin{definition}
\label{42b}
Denote the ``last'' successive minimum by
\begin{equation}
\label{last}
\Lambda(L,\|\;\|) =\lambda_b(L,\|\;\|).
\end{equation}
\end{definition}

\begin{definition}
\label{43a}
A linearly independent family
\[
\{ v_i \}_{i=1,\ldots,b} \subset L
\]
is called {\em quasiorthogonal\/} if~$\|v_i\|= \lambda_i$ for all
$i=1,\ldots,b$.
\end{definition}

Note that a quasiorthogonal family spans a lattice of finite index
in~$L$, but may in general not be an integral basis, a source of some
of the complications of the successive minimum literature.

Dually, we have the Banach space~$B^*=\Hom(B,\R)$, with norm
$\|\;\|^*$, and dual lattice~$L^*\subset B^*$, with~$\rk(L^*)=\rk(L)$.

\begin{theorem}[W. Banaszczyk]
Every lattice~$L$ in every Banach space~$(B,\|\;\|)$ satisfies the
inequality
\begin{equation}
\label{43}
\lambda_1(L,\|\;\|)\, \Lambda(L^*,\|\;\|^*)\le C b(1+\log b),
\end{equation}
for a suitable numerical constant~$C$, where $b=\rk(L)$.
\end{theorem}

In fact, the upper bound is valid more generally for the 
product 
\begin{equation}
\label{szcz}
\lambda_i(L) \lambda_{b-i+1}(L^*),
\end{equation}
for all~$i=1,\ldots,b$ \cite{Bn2}.

\begin{remark}
A lattice~$L\subset \R^b$ admits an {\em orthogonal\/} basis if and
only if~$\lambda_i(L) \lambda_{b-i+1}(L^*) =1$ for all~$i$.  Thus,
Banaszczyk's bound can be thought of as a measure of the
quasiorthogonality of a lattice in Banach space.
\end{remark}

Given a class~$\alpha\in H_k(M;\Z)$ of infinite order, we define the
stable norm~$\| \alpha_\R \|$ by setting
$$
\|\alpha_\R\|=\lim_{m\to \infty} m^{-1} \inf_{\alpha(m)}
\vol_k(\alpha(m)),
$$
where~$\alpha_\R$ denotes the image of~$\alpha$ in real homology,
while~$\alpha(m)$ runs over all Lipschitz cycles with integral
coefficients representing the multiple class~$m\alpha$.  The stable
norm is dual to the comass norm~$\|\;\|^*$ in cohomology,
\cf~\cite{Fe2,BK1}.

\begin{definition}
\label{stsys}
The {\em stable homology~$k$-systole} of~$(X,\gmetric)$ is
\begin{equation}
\label{34}
\stsys_k(\gmetric)=\lambda_1(H_k(X,\Z)_\R,\|\;\|),
\end{equation}
where~$\|\;\|$ is the stable norm.
\end{definition}

\section{Linearity \vs indeterminacy of triple Massey products}
\label{sect:proofs}

We will denote by~$\integerlat{k}$ the image of integral cohomology in
real cohomology under inclusion of coefficients.  Let~$\{[v_i]\}
\subset \integerlat{\two}$ be a quasiorthogonal family in the sense of
Definition~\ref{43a}, with
\[\|v_i\|^*= \lambda_i(\integerlat{m}, \|\;\|^*),
\]
as in formula \eqref{success}, where~$\|\;\|^*$ is the comass norm.
Here we assume, to simplify the calculations, that each~$\two$-form
$v_i$ minimizes the comass norm in its cohomology class.  Given an
exact $(2m)$-form $v_i \wedge v_j$, let $w_{ij}$ be a primitive of
least comass, \cf \eqref{isop}.

\begin{definition}
An element of the form
\[
\left[ w_{ij} \wedge v_k -(-1)^m v_i \wedge w_{jk} \right] \in \langle
v_i, v_j, v_k \rangle
\]
is called an {\em quasiorthogonal element of the Massey
product}~$\langle v_i, v_j, v_k \rangle$.
\end{definition}

\begin{lemma}
\label{51}
Under the hypotheses of Theorem \ref{22}, the existence of a
nontrivial Massey product implies the existence of a nonzero
quasiorthogonal element of a suitable Massey product.
\end{lemma}

\begin{proof}
The lemma follows by linearity, \cf \eqref{45}.  Since the detailed
proof contains a delicate point involving indeterminacy, we include it
here.  By triviality of cup product hypothesis~\eqref{h1a} of
Theorem~\ref{22}, for each pair of indices~$1\le i,j\le b_m(X)$, there
is a~$(2m-1)$-form~$w_{ij}$ solving the equation
\begin{equation}
\label{primitive}
v_i\wedge v_j =dw_{ij}.
\end{equation}
Furthermore, given a metric~$\gmetric$, we can assume that~$w_{ij}$
satisfies the inequality
\begin{equation}
\label{42}
\|w_{ij}\|^* \leq \IQ_\three(\gmetric) \| v_i \wedge v_j \|^*,
\end{equation}
\cf formula \eqref{isop}.

Using index notation (Einstein summation convention), let~$i,j,k$ run
from 1 to~$b_\two (X)$.  Let~$\langle u,v,w \rangle$ be a nontrivial
Massey product, as in Theorem~\ref{22}.  Choose representative
differential forms~$\alpha=\alpha^i v_i \in u$,~$\beta= \beta^j v_j
\in v$, and~$\gamma= \gamma^k v_k \in w$.  Then
\begin{equation}
\begin{split}
\alpha\wedge \beta &= \left( \alpha^i v_i\right) \wedge \left( \beta^j
v_j\right) \\&= \alpha^i \beta^j v_i \wedge v_j \\&= \alpha^i \beta^j
d w_{ij} \\&= d\left( \alpha^i \beta^j w_{ij}\right),
\end{split}
\end{equation}
and similarly~$\beta\wedge \gamma = d\left( \beta^j \gamma^k
w_{jk}\right)$.  Since the Massey product is nontrivial, we obtain a
nonzero cohomology class
\begin{equation}
\label{abc}
\left[\alpha^i \beta^j w_{ij} \wedge \gamma -(-1)^m \alpha \wedge
\beta^j \gamma^k w_{jk} \right] \not =0 \in H_{\dR}^\five(X).
\end{equation}
By linearity, we have
\begin{equation}
\label{45}
\begin{aligned}
\alpha^i \beta^j w_{ij} \wedge \gamma - & (-1)^m \alpha \wedge \beta^j
\gamma^k w_{jk} = \\ & = \alpha^i \beta^j \gamma^k \left(
w_{ij}^{\phantom{I}} \wedge v_k^{\phantom{I}} -(-1)^m_{\phantom{I}}
v_i^{\phantom{I}} \wedge w_{jk}^{\phantom{I}} \right) .
\end{aligned}
\end{equation}
Therefore

\begin{equation}
\alpha^i \beta^j \gamma^k \left[ w_{ij} \wedge v_k -(-1)^m v_i \wedge
w_{jk} \right] \neq 0 \in H_{\dR}^{\five}(X).
\end{equation}
In fact, the nontriviality of the Massey product yields the stronger
conclusion that we have a nonzero class in the quotient
\begin{equation}
\alpha^i \beta^j \gamma^k \left[ w_{ij} \wedge v_k -(-1)^m v_i \wedge
w_{jk} \right] \neq 0 \in H_{\dR}^{\five}(X)/\Indet,
\end{equation}
\cf \eqref{21}.  Hence, for suitable indices~$1\leq s,t,r\leq b_m(X)$,
we obtain a nonzero class
\[
\left[ w_{st} \wedge v_r -(-1)^m v_s \wedge w_{tr} \right] \in 
\langle v_s, v_t, v_r \rangle
\]
in~$H_{\dR}^{\five}(X)/\Indet$.  Note that this conclusion differs
from the assertion that the Massey product $\langle v_s, v_t, v_r
\rangle$ is nontrivial, since its indeterminacy subspace may be
different from that of the Massey product~$\langle u,v,w \rangle$.
\end{proof}

\begin{remark}
The indices $s,t,r$ above may depend on the various choices involved
in the construction, but the key estimate \eqref{62} remains valid,
due to the uniqueness of the least natural number, by the well-ordered
property of $\N$.
\end{remark}

\begin{lemma}
\label{51b}
Let~$\x\in H_{3m-1}(X,\R)$ be a fixed nonzero class.  The hypotheses
of Theorem \ref{22} imply the existence of a nonzero quasiorthogonal
element of a Massey product, which pairs nontrivially with~$\x$.
\end{lemma}

\begin{proof}
Consider the family of all quasiorthogonal elements $q_i$ of Massey
products.  Let~$V$ be the vector space spanned by all such elements
$q_i$.  By~\eqref{45}, the space~$V$ meets every nontrivial Massey
product.  By our Massey-type hypothesis, we have
\begin{equation}
\label{58}
V=H^{3m-1}_{\dR}(X).
\end{equation}

Choose any cohomology class~$a$ which pairs nontrivially with~$\x$,
\ie $a(\x)\not=0$.  By \eqref{58}, we can write $a= a^i q_i$, where
$q_i$ are quasiorthogonal elements of Massey products.  Thus $a^i
q_i(\x) \not=0$ and by linearity, one of the quasiorthogonal elements,
say $q_{i_0}$, also pairs nontrivially with~$\x$.
\end{proof}

\begin{lemma}
Assume~$H^{\four}(X,\Z)$ is torsionfree.  Then every quasior\-thogonal
element of a Massey product satisfies the integrality condition
\begin{equation} 
\label{55}
\int_{\x} \langle v_s, v_t, v_r \rangle \in \Z ,
\end{equation}
where~$x_0\in H_m(X,\Z)$ is any integral class.
\end{lemma}

\begin{proof}
Choose representatives for the~$v_i$ in the cohomology group with
integer coefficients~$H^{\two}(X,\Z)$ in the sense of singular
cohomology theory.  We denote these representatives~$\tilde v_i$.
Choose an~$m$-cocycle~$\tilde {\tilde v}_i \in \tilde v_i$.  Note that
the class
\[
[\tilde {\tilde v}_s\wedge \tilde {\tilde v}_t]\in H^{\four}(X,\Z)
\]
vanishes integrally, and thus the Massey product~$\langle \tilde v_s,
\tilde v_t, \tilde v_r \rangle$ is defined over~$\Z$.  The lemma now
follows from the compatibility of the de Rham and the singular Massey
product theories, verified in \cite{M} and \cite[Section~11]{KR1}, in
terms of differential graded associative (dga) algebras, \cf
Remark~\ref{56} below.
\end{proof}

\begin{remark}
\label{56}
The following three remarks were kindly provided by R.~Hain (see
\cite{KR1, SGT} for more details).

1. If~$A^*$ and~$B^*$ are dga algebras (not necessarily commutative)
and~$f : A^* \to B^*$ is a dga homomorphism that induces an
isomorphism on homology, then Massey products in~$H^*(A^*)$ and
$H^*(B^*)$ correspond under~$f^* : H^*(A^*) \to H^*(B^*)$.

2. If~$M$ is a manifold, then there is a dga~$K^*$ that contains both
the de Rham complex~$A^*(M)$ of~$M$, and also the singular cochain
complex~$S^*(M)$ of~$M$.  The two inclusions
\[
        A^*(M) \to K^* \leftarrow S^*(M)
\]
are both dga quasi-isomorphisms (\ie induce isomorphism in
cohomology), \cf \cite[Corollary~10.10]{FHT}.

3. The point is that the inclusions~$A^*(M) \to K^* \leftarrow S^*(M)$
are both dga homomorphisms (and quasi-isomorphisms), even
though~$A^*(M)$ is commutative and~$S^*(M)$ is not.  Combining these
two remarks, we see that Massey products in singular cohomology and in
de Rham cohomology correspond.  The complex~$K^*$ is a standard tool
in rational homotopy theory.  It is defined as follows.  Let Simp be
the simplicial set of smooth singular simplices of~$M$.  Then~$K^*$ is
Thom-Whitney complex of differential forms on~Simp.
\end{remark}

\section{Proofs of main results}
\label{five}

\begin{proof}[Proof of Theorem~\ref{22}]
Let~$\gmetric$ be a metric on~$X$. Let~$\|\;\|$ be the associated
stable norm in homology.  Choose a class~$\x\in H_\five(X,\Z)_\R$
satisfying
\begin{equation}
\label{66}
\| \x \|=\stsys_\five(X,\gmetric) =\lambda_1(H_\five(X,\Z)_\R,
\|\;\|).
\end{equation}
We can then choose a cohomology class~$\alpha\in \integerlat{3m-1}$
which pairs nontrivially with the class~$\x$, \ie satisfying~$\alpha(
\x)\not=0$.  We will write this condition suggestively as~$\int_{\x}
\alpha \not= 0$.  A reader familiar with normal currents can interpret
integration in the sense of the minimizing normal current representing
the class~$\x$.  Otherwise, choose a rational Lipschitz~$m$-cycle of
volume~$\epsilon$-close to the value \eqref{66}, and let~$\epsilon$
tend to zero at the end of the calculation below.

By Lemma~\ref{51b}, the class~$\alpha$ can be replaced by a
quasiorthogonal element of a Massey product~$\langle v_s, v_t, v_r
\rangle$, which also pairs nontrivially with $\x$.

Recall that~$\|\;\|^*$ is the comass norm in cohomology.  Changing
orientations if necessary, we obtain from \eqref{55} that
\begin{equation}
\label{62}
1 \leq \int_{\x} w_{st} \wedge v_r -(-1)^m v_s \wedge w_{tr},
\end{equation}
and therefore
\begin{equation}
1 \leq C(m) \left( \| w_{st} \|^* \|v_r \|^* + \| v_s \|^* \| w_{tr}
\|^* \right) \| \x \|,
\end{equation}
where~$C(m)$ depends only on~$\two$.  Now by \eqref{42}, we have
\[
\begin{aligned}
1 & \leq 2 C(m) \| v_s\|^* \| v_t \|^* \| v_r\|^*
\IQ_{\three}(\gmetric) \| \x\| \\ & = 2 C(m) \lambda_s \lambda_t
\lambda_r \IQ_{\three}(\gmetric) \| \x\| \\ &\leq 2C(m) \left( \Lambda
\left( \integerlat{\two^{\phantom{i}}}, \|\;\|^* \right) \right)^3
\IQ_{\three}(\gmetric) \| \x\|,
\end{aligned}
\]
by Definition \ref{42b} of the ``last'' successive
minimum~$\Lambda(L)$.  Finally, by definition we
have~$\stsys_\two(\gmetric) = \lambda_1(H_\two(X), \|\;\|)$,
where~$\|\;\|$ is the stable norm, and therefore
\begin{equation}
\stsys_\two(\gmetric)^3 \leq 2C(m) \left( \lambda_1^{\phantom{A}}
(H_\two(X))\, \Lambda(H^\two(X)) \right)^3 \IQ_{\three}(\gmetric) \|
\x\|.
\end{equation}
Applying Banaszczyk's inequality \eqref{43}, we obtain
\[
\begin{aligned}
\stsys_\two(\gmetric)^3 & \leq C(m) (b(1+\log b))^3
\IQ_{\three}(\gmetric) \| \x\| \\ & = C(m) (b(1+\log b))^3
\IQ_{\three}(\gmetric) \stsys_{\five}(\gmetric),
\end{aligned}
\]
where~$b=b_{\two}(X)$, while the new coefficient~$C(m)$ incorporates
the numerical constant from Banaszczyk's inequality.  This completes
the proof of Theorem~\ref{22}.
\end{proof}

\begin{proof}[Proof of Theorem \ref{222}]
Exploiting the orientability of~$X$, we represent its fundamental
cohomology class as a product~$\langle u_1, u_2, u_3 \rangle \cup
u_4$, with~$u_i \in H^2(X)$.  Here we write~$\langle u_1, u_2, u_3
\rangle$ as shorthand for an orthogonal element of a Massey product,
while~$u_4$ may be chosen to be any class which pairs nontrivially
with the Poincar\'e dual of~$\langle u_1, u_2, u_3 \rangle$.
Relation~\eqref{55} is replaced by the following integrality relation
among the elements~$v_i\in H_{\dR}^2(X)$ of a quasiorthogonal family:
\begin{equation} 
\label{65}
\int_X \langle v_s, v_t, v_r \rangle \cup v_p \in \Z \setminus \{ 0
\}.
\end{equation}
The rest of the proof is similar.
\end{proof}

\begin{proof}[Proof of Proposition~\ref{81}]
Choose a class~$\x\in H_8(X,\Z)_\R$ satisfying $\| \x \| = \lambda_1(
H_8(X,\Z)_\R; \|\;\|)$.  The class~$\x$ pairs nontrivially with one of
the classes~$u\cup v$ or~$w$.

If for some Massey product~$u$, we have~$\int _{\x} u \cup v\not=0$,
we argue as in the proof of Theorem~\ref{22}, exploiting the
hypothesis that the cup product in~$H^2(X)$ is trivial, in order to
define the quasiorthogonal elements of triple Massey products.

If the class~$w$ satisfies~$w( \x)\not=0$, we argue with a
quasiorthogonal family in~$\integerlat{4}$ as in~\cite{BK1} to obtain
the lower bound for the stable norm of~$ \x$ in terms
of~$\stsys_4(\gmetric)^2$.
\end{proof}

\section{Acknowledgments}  

The author is grateful to R.~Hain for clarifying the compatibility of
the singular and the de Rham theories of Massey products, \cf
Remark~\ref{56}; to Y. Rudyak for carefully reading an earlier version
of the manuscript and making valuable comments; to A.~Suciu for
helpful discussions and for developing examples of manifolds with
nontrivial Massey triple products; to Sh.~Weinberger for stimulating
conversations; and to the referee for catching an error in an earlier
version of the text.

\bibliographystyle{amsalpha}

\begin{thebibliography}{ABCDE}

\bibitem[Am04]{Am} Ammann, B.: Dirac eigenvalue estimates on two-tori.
{\em J. Geom. Phys.} \textbf{51} (2004), no. 3, 372--386.

\bibitem[BabB05]{BB} Babenko, I.; Balacheff, F.: G\'eom\'etrie
systolique des sommes connexes et des rev\^etements cycliques, {\em
Math. Annalen\/} \textbf{333} (2005), no.~1, 157-180.

\bibitem[Ban93]{Bn1} Banaszczyk, W.: New bounds in some transference
theorems in the geometry of numbers, {\em Math. Ann.} \textbf{296}
(1993), 625--635.

\bibitem[Ban96]{Bn2} Banaszczyk, W.: Inequalities for convex bodies
and polar reciprocal lattices in~$\R^n$. \textup{II}. Application of
$K$-convexity, {\em Discrete Comput. Geom.} \textbf{16} (1996),
305--311.


\bibitem[BCIK05]{BCIK1} Bangert, V; Croke, C.; Ivanov, S.; Katz, M.:
Filling area conjecture and ovalless real hyperelliptic surfaces, {\em
Geometric and Functional Analysis (GAFA)\/} \textbf{15} (2005), no.~3,
577-597.  See \texttt{arXiv:math.DG/0405583}

\bibitem[BCIK07]{BCIK2} Bangert, V; Croke, C.; Ivanov, S.; Katz, M.:
Boundary case of equality in optimal Loewner-type inequalities.  {\em
Trans. Amer. Math. Soc.} \textbf{359} (2007), no.~1, 1--17.  See
\texttt{arXiv:math.DG/0406008}



\bibitem[BK03]{BK1} Bangert, V.; Katz, M.: Stable systolic
inequalities and cohomology products, {\em Comm. Pure Appl. Math.}
\textbf{56} (2003), 979--997.  Available at
\texttt{arXiv:math.DG/0204181}

\bibitem[BK04]{BK2} Bangert, V; Katz, M.: An optimal Loewner-type
systolic inequality and harmonic one-forms of constant norm.  {\em
Comm. Anal. Geom.}  \textbf{12} (2004), no. 3, 703-732.  See
\texttt{arXiv:math.DG/0304494}


\bibitem[BKSS06]{e7} Bangert, V; Katz, M.; Shnider, S.; Weinberger,
S.: $E_7$, Wirtinger inequalities, Cayley 4-form, and homotopy.  See
\texttt{arXiv:math.DG/0608006}




\bibitem[CLOT03]{CLOT} Cornea, O.; Lupton, G.; Oprea, J.; Tanr\'e, D.:
Lusternik-Schnirelmann category.  {\em Mathematical Surveys and
Monographs}, \textbf{103}.  Amer. Math. Soc., Providence, RI, 2003.

\bibitem[CrK03]{CK} Croke, C.; Katz, M.: Universal volume bounds in
Riemannian mani\-folds, {\it Surveys in Differential Geometry\/}
\textbf{8} (2003), 109-137.  Available at
\texttt{arXiv:math.DG/0302248}


\bibitem[DR03]{DR} Dranishnikov, A.; Rudyak, Y.: Examples of
non-formal closed~$(k-1)$-connected manifolds of dimensions~$\ge
4k-1$, {\em Proc. Amer. Math. Soc.} \textbf{133} (2005), no. 5,
1557-1561.  See \texttt{arXiv:math.AT/0306299}.

\bibitem[Fe74]{Fe2} Federer, H.: Real flat chains, cochains, and
variational problems, {\em Indiana Univ. Math. J.} \textbf{24} (1974),
351--407.


\bibitem[FHT98]{FHT} Felix, Y.; Halperin, S.; Thomas, J.-C.: Rational
homotopy theory.  {\em Graduate Texts in Mathematics}, \textbf{205}.
Springer-Verlag, New York, 2001.


\bibitem[Gr83]{Gr1} Gromov, M.: Filling Riemannian manifolds.  {\em
J. Differential Geom.}  \textbf{18} (1983), 1--147.


\bibitem[Gr96]{Gr2} Gromov, M.: Systoles and intersystolic
inequalities, Actes de la Table Ronde de G\'{e}om\'{e}trie
Diff\'{e}rentielle (Luminy, 1992), 291--362, {\em S\'{e}min. Congr.},
\textbf{1}, Soc. Math. France, Paris, 1996.
\newline\noindent
\texttt{www.emis.de/journals/SC/1996/1/ps/smf\_sem-cong\_1\_291-362.ps.gz}


\bibitem[Gr99]{Gr3} Gromov, M.: Metric structures for Riemannian and
non-Riemannian spaces, {\em Progr. in Mathematics\/} \textbf{152},
Birkh\"{a}user, Boston, 1999.

\bibitem[He86]{He} Hebda, J.: The collars of a Riemannian manifold and
stable isosystolic inequalities, {\em Pacific J. Math.} \textbf{121}
(1986), 339--356.

\bibitem[IK04]{IK} Ivanov, S.; Katz, M.: Generalized degree and
optimal Loewner-type inequalities, {\em Israel J. Math.}  \textbf{141}
(2004), 221-233.  \texttt{arXiv:math.DG/0405019}

\bibitem[Ka03]{Ka3} Katz, M.: Four-manifold systoles and surjectivity
of period map, {\em Comment. Math. Helv.} \textbf{78} (2003), 772-786.
See \texttt{arXiv:math.DG/0302306}


\bibitem[Ka07]{SGT} Katz, M.: Systolic geometry and topology.  {\em
Mathematical Surveys and Monographs}, to appear.  American
Mathematical Society, Providence, R.I.




\bibitem[KL05]{KL} Katz, M.; Lescop, C.: Filling area conjecture,
optimal systolic inequalities, and the fiber class in abelian covers.
Geometry, spectral theory, groups, and dynamics, 181--200, {\em
Contemp. Math.} \textbf{387}, Amer. Math. Soc., Providence, RI, 2005.
See \texttt{arXiv:math.DG/0412011}


\bibitem[KR06]{KR1} Katz, M.; Rudyak, Y.: Lusternik-Schnirelmann
category and systolic category of low dimensional manifolds.  {\em
Communications on Pure and Applied Mathematics\/} \textbf{59} (2006),
no.~10, 1433-1456.  Available at \texttt{arXiv:math.DG/0410456}


\bibitem[KR07]{KR2} Katz, M.; Rudyak, Y.: Bounding volume by systoles
of 3-manifolds.  See \texttt{arXiv:math.DG/0504008}



\bibitem[KRS06]{KRS} Katz, M.; Rudyak, Y.; Sabourau, S.: Systoles of
2-complexes, Reeb graph, and Grushko decomposition.  {\em
International Math. Research Notices\/} \textbf{2006} (2006).  See
\texttt{arXiv:math.DG/0602009}



\bibitem[KS05]{KS2} Katz, M.; Sabourau, S.: Entropy of systolically
extremal surfaces and asymptotic bounds, {\em Ergodic Theory and
Dynamical Systems}, \textbf{25} (2005), no.~4, 1209-1220.  Available
at \texttt{arXiv:math.DG/0410312}

\bibitem[KS06a]{KS1} Katz, M.; Sabourau, S.: Hyperelliptic surfaces
are Loewner, {\em Proc. Amer. Math. Soc.} \textbf{134} (2006), no.~4,
1189-1195.  Available at the site \texttt{arXiv:math.DG/0407009}


\bibitem[KS06b]{KS3} Katz, M.; Sabourau, S.: An optimal systolic
inequality for CAT(0) metrics in genus two.  {\em Pacific J. Math.}
\textbf{227} (2006), no.~1, 95-107.  See
\texttt{arXiv:math.DG/0501017}


\bibitem[KSV05]{KSV} Katz, M.; Schaps, M.; Vishne, U.: Logarithmic
growth of systole of arithmetic Riemann surfaces along congruence
subgroups.  {\em J. Differential Geom.} (2007).  Available at
\texttt{arXiv:math.DG/0505007}



\bibitem[Ma69]{M} May, J.: Matric Massey products.  {\em J. Algebra\/}
\textbf{12} (1969) 533--568.

\bibitem[RT00]{RT} Rudyak, Y.; Tralle, A.: On Thom spaces, Massey
products, and nonformal symplectic manifolds.  {\em
Internat. Math. Res.  Notices\/} \textbf{10} (2000), 495--513.

\bibitem[Sa04]{Sa1} Sabourau, S.: Systoles des surfaces plates
singuli\`eres de genre deux, {\em Math. Zeitschrift\/} \textbf{247}
(2004), no. 4, 693--709.




\bibitem[Sa06]{Sa3} Sabourau, S.: Systolic volume and minimal entropy
of aspherical manifolds, {\em J. Differential Geom.}, \textbf{74}
(2006), no.~1, 155-176.  Available at \texttt{arXiv:math.DG/0603695}



\bibitem[Si05]{Sik} Sikorav, J.: Bounds on primitives of differential
forms and cofilling inequalities.  See \texttt{arXiv:math.DG/0501089}

\end{thebibliography}

\end{document}